\documentclass[11pt]{amsart}
\usepackage{geometry}                
\usepackage{amssymb}
\usepackage{enumerate}
\usepackage{color}
\usepackage{graphics}
\usepackage{epsfig}

\newcommand{\real}{\mathbb{R}}

\newtheorem{theorem}{Theorem}[section]

\theoremstyle{definition}
\newtheorem{definition}[theorem]{Definition}
\newtheorem{example}[theorem]{Example}

\newtheorem{proposition}[theorem]{Proposition}

\theoremstyle{remark}

\numberwithin{equation}{section}
\begin{document}

\title{On $p$-Parabolicity of Riemannian Submersions}
\author{Maria Andrade}
\email{maria.ufal@gmail.com}

\author{Pietro da Silva}
\email{pietro.ufal@gmail.com}


\begin{abstract}
We provide some criteria to $p$-parabolicity of Riemannian submersions. In particular, if $N$ is $p$-parabolic and $\pi:M\to N$ is a Riemannian submersion with uniformly bounded volume of fibers, then $M$ is also $p$-parabolic. In the  case of warped manifolds we characterize $p$-parabolicity in terms of a volume growth condition.
\end{abstract}

\maketitle
\section{Introduction}

Let $(M,g)$ be a connected Riemannian manifold, and a pair of subsets $D\subset\Omega\subset M$ with $D$ compact and $\Omega$ a connected domain. Given  $p\in (1,\infty),$ the $p$-capacity of $D$ in $\Omega$ is defined by
 \begin{equation}\label{pcapacity}
 \text{Cap}_p(D,\Omega):=\text{inf}\left\{\int_{\Omega} |\nabla u|^p: u\in C_0^1(\Omega),\ u\geq 1\  \text{on}\ D \ \right\}.
 \end{equation}
If $\Omega=M$ then we write $\text{Cap}_p(D)$ for simplicity. Due to well known properties of $p$-capacity, if $\Omega_t\subset\Omega_{t+1}$ and $\cup_{t=1}^{\infty}\Omega_t=M$, then $\text{Cap}_p(D)=\lim_{t\to\infty}\text{Cap}_p(D,\Omega_t)$. We say that $M$ is $p$-\textit{parabolic} if $\text{Cap}_p(D)=0$ for any compact $D\subset M$ and $p$-\textit{hyperbolic} otherwise.

The $p$-parabolicity is closely related with properties of the $p$-Laplacian operator $\triangle_pu=\text{div}(|\nabla f|^{p-2}\nabla u)$, defined over real functions on $M$. In fact, the $p$-Laplacian is the Euler-Lagrange operator associated to the energy functional in the  right side of (\ref{pcapacity}) (see \cite{Troyanov00}, \cite{CHS01}). The case $p=2$  has been extensively studied linking several mathematical areas, namely geometry, analysis and probability (\cite{Grigoryan99} provides a deep survey on this topic). The classical Laplace-Beltrami operator $\triangle=\triangle_2$ carries many of geometric aspects of $(M,g)$. On the other hand, equations involving $\triangle$ are subject of several analytical problems including the heat equation whose solution can be used to generate the Brownian motion on $M$. 2-parabolicity is equivalent to non-existence of a non-constant superharmonic ($\triangle u\leq0$) function on $M$ as well as recurrence of the Brownian Motion on $M$.

In this paper we deal with Riemannian submersions $\pi:M\to N$ and we attempt to know conditions on $N$ to assure $p$-parabolicity on $M$. More precisely, our first result reads as follows:

\begin{theorem} \label{volim}
 Let $M$ and $N$ be Riemannian manifolds and $\pi:M\to N$ a Riemannian submersion with fibers $\mathcal{F}_x=\pi^{-1}(x)$ having uniformly bounded volume, i.e, $Vol(\mathcal{F}_x)\leq C$. If $N$ is $p-$parabolic, then M is $p-parabolic.$ 
\end{theorem}

A related work was made by \cite{BrandaoOliveira13} that asserts the 2-parabolicity of a manifold submersed in a parabolic one in such a way that fibers are minimal and compact. In addition to dealing with the more general setting of $p$-parabolicity, here we drop conditions over mean curvature to use weaker restrictions on volume of fibers. 

In the special case of warped products we are able to claim necessary and sufficient conditions to $p$-parabolicity.
\begin{theorem}\label{main}
Consider the warped product $M=N\times_{f}L$, where $N^n$ is Riemannian manifold admitting exhaustion by closed geodesic balls $(B_t)_{t>0}$ centered in the same point, $L^{\ell}$ has bounded volume and $p>1$. A necessary and sufficient condition to $M$ be $p$-parabolic is
\begin{equation}\label{hypothesis}
\int_1^{\infty}\left(\int_{\partial B_t}f(x)^{\ell}\mu'_N(dx)\right)^{\frac{1}{1-p}}dt=\infty. 
\end{equation}
\end{theorem}

Our necessary and sufficient condition to $p$-parabolicity of warped products generalizes \cite{TroyanovSAM99} which treats cylindrical warped manifolds. In the Theorem~\ref{main} we can consider $M=D\cup (N\times_fL)$, where $D$ is precompact subset of $M$.

This paper is organized in this way: section 2 states basics definitions to be used in the following; section 3 presents examples to situate our results in the literature; proofs are in section 4.

\section{Prelimiaries}

\begin{definition}A smooth map $\pi: (M,g)\to (N,h)$ is a Riemannian submersion if $\pi_*$ is surjective and satisfies the following property: $$g_x(v,w)=h_{\pi(x)}(\pi_*v,\pi_*w)$$ for any $v,w$ tangent vectors in $T_xM$ and perpendicular to the kernel of $\pi_*$.
\end{definition}

%

Note that if $N$ and $L$ be Riemannian manifolds, then $\pi:N\times L\to L$ is a Riemannian submersion, because, for each point $l\in L,$ the application $\pi|_{N\times l}$ is an isometry under N. Moreover, an important example of Riemannian submersion is called warped product, introduced by \cite{Bishop69}, which we define bellow.

\begin{definition} (Warped Product) Let $N$ and $L$ be Riemannian manifolds and $f:N\to(0,\infty)$ a differentiable function. Consider the product (differentiable) manifold $N\times L$ with its projections $\pi_1: N\times L \to N$ and $\pi_2:N\times L \to L.$ The {\it warped product} $W=N\times_f L$ is the manifold
$N\times L$ furnished with the Riemannian structure such that
\begin{eqnarray}
\label{metricawp}
ds^2_W=ds^2_N+f^2(n)ds^2_L,
\end{eqnarray}
where $n\in N.$

\vspace{0.1cm}

We call $N$ and $L$ basis (or leaves) and fibers, respectively.
\end{definition}
%
%

\section{Examples}

In this section we show some examples to situate the range of our results in the literature. Before, remember that $\real$ is $1-$parabolic and $\real^n$ is $p-$parabolic for $p\geq n$ as consequence of its decomposition in $B_1\cup([1,\infty)\times_f\mathbb{S}^{n-1})$, where $B_1$ is the closed unitary ball in $\real^n$ and $f(x)=x$. In the following examples $N=\real^n$ and $M=N\times_f L$.

\begin{example}
Let $\pi:M\to N,$ be a Riemannian submersion with $L= \mathbb{S}^n, \ f= e^{-x^2},$ then by the Theorem \ref{volim} $\real\times_f \mathbb{S}^n$ is $2-$parabolic. Note that our result generalizes the Theorem 2.3 in \cite{BrandaoOliveira13} since in this example the fibers are not minimals. 
\end{example}

Differently than \cite{BrandaoOliveira13} asserts, their compactness assumption plays no essential rule to assure parabolicity of $M$. In fact, a hypothesis over volume of fibers is what matters as explained in the next example.

\begin{example}
 Consider $L$ a non-compact finite volume Riemannian manifold. If $f$ is a bounded function, then by Theorem~\ref{main} $M$ is $p-$parabolic for $p\geq n$. exposing  general Theorem~\ref{volim} is to relate $p-$parabolicity over submersions, this example can be applied 
\end{example}

Theorem~\ref{main} assures that, satisfied the volume condition on $L$, the only geometrical characteristic of $L$ that interferes in the $p-$parabolicity of $M$ is its dimension, delivering to $N$ and $f$

\section{Proofs}
To prove theorem~\ref{volim} we need the following result:

\begin{proposition}[Troyanov]
\label{troyanov_sequence}
The domain $\Omega$ is $p-$parabolic if and only if there exists a sequence of functions $u_j \in C_0^1(\Omega)$ such that $0\leq u_1\leq 1,$ $u_j\to 1$ uniformly on every compact subsets of $\Omega$ and 
 \begin{eqnarray*}
  \int_{\Omega}|\nabla u_j|^p\to 0.
 \end{eqnarray*}
\end{proposition}

\proof{ (Theorem~\ref{volim}) Since N is $p-$ parabolic by Proposition \ref{troyanov_sequence} there exists a sequence of functions $u_j\in C_0^1(N)$ such that $0\leq u_j\leq 1,\ u_j\to 1$ uniformly on every compact of $N$ and $$\int_{N}|du_j|^p\to 0.$$

Consider $\widetilde{u}_j=u_j\circ \pi.$ Note that, $\{\widetilde{u}_j\}\subset C_0^1({M}),\ 0\leq \widetilde{u_j}\leq 1,$ and $\widetilde{u_j}\to 1$ uniformly on every compact of $M$, because if $K\subset{M}$ is compact and  $y\in K$, then $\pi(y)\in\pi(K)$, a compact subset of $N$, so $|\widetilde{u}_j(y)-1|=|u_j(\pi(y))-1|\to 0.$ Moreover,
\begin{eqnarray*}
 \int|d\widetilde{u}_j|^p&=&\int_{M}|\nabla\widetilde{u}_j|^pd\mu_M\\
 &=&\int_{N}\int_{\mathcal{F}_x}|\nabla{u}_j|^p d\mathcal{F}_xd\mu_N\\
&=& \int_{N}Vol(\mathcal{F}_x)|\nabla\widetilde{u}_j|d\mu_N\\
&\leq& C\int_{N}|du_j|^pd\mu_N\xrightarrow{j\to{\infty}}0.
\end{eqnarray*}
So, by proposition 4.1 in \cite{TroyanovSAM99} $M$ is $p-$parabolic.
}
\endproof

To prove the theorem~\ref{main} we need the definition of $p-$flux of a function $h$.  

\begin{definition}To each pair $D\subset\subset\Omega\subset M$ the class $\Lambda(D,\Omega)$ of functions $h:\Omega \to \real$ such that
\begin{enumerate}[i)]
 \item $h$ is continous, locally Lipschitz, non constant and bounded below;
 \item $D\subset\{x\in\Omega:h(x)=r_0:=min h\};$
 \item if $r<\text{sup }h$ then $\{x\in\Omega:h(x)\leq r\}$ is compact.
\end{enumerate}

The $p-$flux of a function $h\in \Lambda(D,\Omega)$ is the function $\Phi_{h,p}: [r_0,r_1)\to\real$ defined by
\begin{eqnarray*}
 \Phi_{h,p}=\int_{\partial \Omega_r}|\nabla h(x)|^{p-1}\mu'(dx),
\end{eqnarray*}

where $\Omega_r:=\{x\in\Omega:h(x)<r\}$ is the range of $h$ (i.e., $r_0:=\text{min} h \in \real$) and $r_1:=\text{sup} h \in \real\cup\{\infty\}.$
\end{definition}

We use the next result to prove the theorem~\ref{main} that shows relation between $p-$capacity and $p-$flux.

\begin{theorem}[Troyanov]
 Let $D\subset\subset\Omega\subset M$ and $p>1,$ then
 \begin{eqnarray*}
  \text{Cap}_p(D,\Omega)=inf_{h\in\Lambda(D,\Omega)}\left(\int_{r_0}^{r_1}\Phi_{h,p}(r)^\frac{1}{1-p}dr\right)^{1-p}.
 \end{eqnarray*}
\end{theorem}

\proof{(Theorem~\ref{main}) To prove Theorem~\ref{main} we going to construct subsets $D,D_R\subset W$ such that $M=\cup_{R>1}D_R$ with the $p$-capacity in the form $$Cap_p(D,D_R)=vol(N)\left[\int_1^{R}\left(\int_{\partial B_t}f(x)^n\mu'_n(dx)\right)^{\frac{1}{1-p}}dt\right]^{1-p}.$$ Therefore, the desired equivalence will follow from $$Cap_p(D,W)=\lim_{R\to\infty}Cap_p(D,D_R)$$ (see \cite{TroyanovSAM99} for details).

Take $x_0\in M$ and let $B_r\subset M$ be the geodesic ball with radius $r$ centered in $x_0$. Set $D=B_1\times_fN$ and $D_R=B_R\times_fN$ for $R>1$ and define $h:W\to[1,\infty)$ by taking $h|_D\equiv1$ and $h(z)=\rho(\pi(z))$ if $z\notin D$, where $\rho$ represents the distance to $x_0$. Note that $h\in\Lambda(D,D_R)$ for every $R>1$, and $$|\nabla h(z)|=|\pi_*(z)(\nabla\rho(\pi(z))|=|\nabla\rho(\pi(z))|\leq1.$$ 

Thus, the $p$-flux $\Phi_{h,p}:[1,R)\to \real$ satisfies 

\begin{eqnarray*}
\Phi_{h,p}(r)&=&\int_{\partial D_r}|\nabla h(z)|^{p-1}\mu'_W(dz)\\
&\leq&\int_{\partial D_r}\mu'_W(dz)=\int_{\partial B_r}\left(\int_Nd_{f(x)}\mu_N\right)\mu'_M(dx)\\
&=&\int_{\partial B_r}vol(N)f(x)^n\mu'_M(dx).
\end{eqnarray*}

 Therefore,
\begin{eqnarray*}
 Cap_p(D,D_R)&=&\inf_{\tilde{h}\in\Lambda(D,D_R)}\left(\int_1^R\Phi_{\tilde{h},p}(r)^{\frac{1}{1-p}}dr\right)^{1-p}\label{troyanov5.1}\\
 &\leq& \left(\int_1^R\Phi_{h,p}(r)^{\frac{1}{1-p}}dr\right)^{1-p}\nonumber\\
 &\leq& vol(N)\left[\int_1^R\left(\int_{\partial B_r}f(x)^n\mu'_Md(x)\right)^{\frac{1}{1-p}}dr\right]^{1-p}.\nonumber
\end{eqnarray*}

To attain the converse inequality consider a polar coordinate system on $M$ with origin in $x_0$, then a point $x\in M$ will be viewed as a pair $(t,x')$ where $t=\rho(x)$ and $x'\in\partial B_t$, whereas $z\in W$ is a triple $(t,x',y)$. Let us take an arbitrary test function $u\in C_0^1(D_R)$ such that $u|_D\equiv1$. Then we have
\begin{eqnarray*}
 1&=&\left|\int_1^R\frac{\partial u(t,x',y)}{\partial t}dt\right|\leq\int_1^R|\nabla u(t,x',y)|dt\\
 &=&\int_1^R|\nabla u(t,x',y)|\left(\int_{\delta B_t}f(x)^n\mu'_M(dx)\right)^{\frac{1}{p}}\left(\int_{\delta B_t}f(x)^n\mu'_M(dx)\right)^{-\frac{1}{p}}dt\\
 &\le&\left[\int_1^R\left(|\nabla u(t,x',y)|^p\int_{\partial B_t}f(x)^n\mu'_M(dx)\right)dt\right]^{\frac{1}{p}}\cdot\\
 &\cdot&\left[\int_1^R\left(\int_{\partial B_t}f(x)^n\mu'_M(dx)\right)^{\frac{1}{1-p}}dt\right]^{\frac{p-1}{p}}
\end{eqnarray*}

in another way $$\left[\int_1^R\left(\int_{\partial B_t}f(x)^n\mu'_M(dx)\right)^{\frac{1}{1-p}}dt\right]^{1-p}\leq\int_1^R\left(|\nabla u(t,x',y)|^p\int_{\partial B_t}f(x)^n\mu'_M(dx)\right)dt$$
and integrating over $N$ we get
\begin{eqnarray*}
 vol(N)I(R,p,t)&\leq&\int_N\int_1^R\left(|\nabla u(t,x',y)|^p\int_{\partial B_t}f(x)^n\mu'_M(dx)\right)dt\mu_N(dy)\\
 &=&\int_{D_R}|\nabla u(z)|^p\mu_W(dz),
\end{eqnarray*}
where $$I(R,p,t)=\left[\int_1^R\left(\int_{\partial B_t}f(x)^n\mu'_M(dx)\right)^{\frac{1}{1-p}}dt\right]^{1-p}.$$

Finally, by taking the infimum among test functions $u$ we obtain the desired inequality.
}
\endproof
\bibliographystyle{plain}
\bibliography{references.bib}

\end{document}